\documentclass[11pt]{amsart}
\usepackage{amssymb,amsmath,amsthm}
\usepackage[hmargin=3cm,vmargin=3.5cm]{geometry}
\usepackage{graphicx}
\usepackage{amscd}
\usepackage[all, knot]{xy}
\usepackage{epsfig}
\usepackage{psfrag}
\newtheorem{theorem}{Theorem}[section]

\newtheorem{definition}[theorem]{Definition}
\newtheorem{con}[theorem]{Conjecture}
\newtheorem{proposition}[theorem]{Proposition}

\title{Quasi-alternating links and odd homology: computations and conjectures}

\author{SLAVIK JABLAN AND RADMILA SAZDANOVI\'C}

\begin{document}
\maketitle

\begin{abstract}
This paper contains computational results about quasi-alternating knots and links and their odd Khovanov homology obtained using by A. Shumakovitch's program \emph{KhoHo}, the program \emph{Knotscape}, and our program \emph{LinKnot}. We enumerate all quasi-alternating links up to 12 crossings. We give the first examples of quasi-alternating knots and links which have at least two different minimal diagrams, one quasi-alternating and the other not quasi-alternating. We also give examples of knots and links with $n\le 12$ crossings which are homologically thin and have no minimal quasi-alternating diagrams. These links are potentially homologically thin links that are not quasi-alternating. J.~Greene  [Gr] proved that one such knot  [JaSa1] $11n_{50}$ is indeed an example.
\end{abstract}


\section{Introduction}

In this paper we present computational results addressing the
classification of quasi-alternating (short QA) links based on the
thickness of Khovanov homology $\overline{Kh}$ \cite{Kh1} and odd
Khovanov homology $\overline{Kh}'$ \cite{OzRaSz} and analyzing their
minimal diagrams.

Our motivation stems from the results obtained from the classical
Khovanov homology $\overline{Kh}$. D.~Bar-Natan was the first to
notice \cite{BN} that the vast majority of prime knots up to 10
crossings (238 among 250) are $\overline{Kh}$-thin. In particular,
all alternating links are $\overline{Kh}$-thin \cite{Le} and all
adequate non-alternating knots are $\overline{Kh}$-thick \cite{Kh2}.
The complete list of $\overline{Kh}$-thick knots up to 13 crossings
is computed by A.~Shumakovitch \cite{Sh1}.

Moreover, C.~Manolescu and P.~Ozsv\'ath show that both Khovanov
homology $\overline{Kh}$ and Heegaard-Floer homology $\widehat{HFK}$
can be used to detect links which are not QA. Quasi-alternating
links are Khovanov homologically $\sigma $-thin (over $\mathbb{Z}$)
and Floer homologically $\sigma $-thin (over
$\mathbb{Z}$/2$\mathbb{Z}$) \cite{MaOz}. The same property extends
to odd Khovanov homology $\overline{Kh}'$ \cite{OzRaSz} and we will
use this property in the rest of the paper.   A knot or link is
called {\it homologically thin} (without qualification), if it is
simultaneously thin with respect to $\overline{Kh}$,
$\widehat{HFK}$, and $\overline{Kh}'$ [Gr, Def. 1.2].

\begin{definition}
The set $Q$ of quasi-alternating links is the smallest set of
links such that

\begin{itemize}
\item the unknot is in $Q$;
\item if the link $L$ has a diagram
$D$ with a crossing $c$ such that

\begin{enumerate}
\item both smoothings of $c$, $L_0$ and $L_{\infty}$, are in $Q$;
\item det($L$)=det$(L_0)$+det$(L_{\infty})$
\end{enumerate}

\end{itemize}

\noindent then $L$ is in $Q$. We say that a crossing $c$
satisfying the properties above is a quasi-alternating crossing of
the diagram $D$ or that $D$ is quasi-alternating at the crossing
$c$ \cite{OzSz,ChKo1}. \label{def}
\end{definition}

The recursive definition makes it difficult to determine if a knot
is quasi-alternating.  It is a challenge to find candidates for
homologically thin knots that are not QA. For a long time, knots
$9_{46}$ = $3,3,-3$ and $10_{140}$ = $4,3,-3$ have been the main
candidates. However, according to A.~Shumakovitch's computations
\cite{Sh2} they are not QA, since are $\overline{Kh}'$-thick,
although they are both $\widehat{HFK}$ and $\overline{Kh}$-thin.

According to Theorem 1 [ChKo1], quasi-alternating links with a
higher number of crossings can be obtained as extension of links
which are already recognized as quasi-alternating \cite{ChKo1,Wi}.

Consider the crossing $c$ in Definition 1 as 2-tangle with marked
endpoints. Using Conway's notation for rational tangles, let
$\varepsilon (c)=\pm 1$, according to whether the overstrand has
positive or negative slope. We will say that a rational 2-tangle
$\tau =C(a_1,\ldots ,a_m)$ extends $c$ if $\tau $ contains $c$ and
$\varepsilon (c)\cdot a_i\ge 1$ for $i=1,\ldots, m$. In particular,
$\tau $ is an alternating rational tangle.

\begin{theorem}
If $L$ is a quasi-alternating link, let $L'$ be obtained by
replacing any quasi-alternating crossing $c$ with an alternating
rational tangle that extends $c$. Then $L'$ is quasi-alternating
$[ChKo1]$.
\end{theorem}

In this paper we give new computational results for QA links up to
12 crossings and the examples of QA links with at least two
different minimal diagrams, where one is QA and the other is not. We
provide examples of knots and links (short $KL$s) with $n\le 12$
crossings which are homologically thin and have no minimal
quasi-alternating diagrams. In the first version of this paper we
proposed these $KL$s as the  candidates for prime homologically thin
links that are not QA, and J.~Greene proved that knot $11n_{50}$ is
not QA [Gr]. Using the method described in his paper it can be shown
that the link $L11n_{90}$ is not QA, although it is homologically
thin. The remaining candidates may require additional ideas to prove
that they are non-QA, if indeed this is the case.  \\

Knots and links are given in Conway notation \cite{Con, JaSa}, which
is implemented in {\it Mathematica} package {\it LinKnot} used for
deriving families of $KL$s and their distinct diagrams. Odd Khovanov
homology $\overline{Kh}'$ is computed using \emph{KhoHo} by A.
Shumakovitch \cite{Sh1}. All flype-equivalent minimal diagrams of
non-alternating $KL$s up to $n \le 12$ crossings are derived from
alternating link diagrams, using software \emph{LinKnot}. {\it
KnotFind}, the part of the program {\it Knotscape} \cite{HosThi}, is
used for recognition of knots, and Jones and Kauffman polynomials
are used for distinguishing links. In addition, we used the
criterion that homologically thick knots are not QA and that for a
quasi-alternating crossing both smoothings must be homologically
thin knots or links.

\section{Quasi-alternating knots up to 12 crossings}

Since an homologically thick knot cannot be QA, we first selected
knots which are homologically thin. Table 1 gives an overview of the
numbers of non-alternating knots with $8 \le n\le 12$ crossings and
how many among them are homologically thin:
\begin{table}[h]
\begin{tabular}{|l|c|c|c|c|c|}   \hline
   No. of crossings & 8 & 9 & 10 & 11 & 12 \\ \hline
   No. of non-alternating knots & 3 & 8 & 42 & 185 & 888 \\ \hline
   No. of homologically thin non-alternating knots & 2 & 6 & 31 & 142 & 663 \\ \hline
\end{tabular}
\vspace{0.2cm}
\caption{ Numbers of non-alternating knots and those non-alternating that are also nomologically thin, listed according to their minimal crossing number.} \label{T1}
\end{table}

Among $11$-crossing knots,  six knots  in  Table \ref{T22} are
$\overline{Kh}$-homologically thin \cite{Sh1} and have a minimal
diagram which is not QA.  Even columns in Table 2 contain the Conway symbols of these knots.
Four of them, $9_{46}$, $10_{140}$, $11n_{139}$, and $11n_{107}$ are
$\overline{Kh}'$-thick \cite{Sh2}, so they are not QA.

\begin{table}[h]
\begin{tabular}{|c|c|c|c|} \hline
 $9_{46}$  & $3,3,-3$ & $10_{140}$  & $4,3,-3$ \\ \hline
 $11n_{139}$  & $5,3,-3$ & $11n_{107}$ & $-2\,1\,2,3,3$ \\ \hline
 $11n_{50}$  & $-2\,2,2\,2,3$ & $11n_{65}$ & $(3,-2\,1)\,(2\,1,2)$ \\ \hline
\end{tabular}
\vspace{0.2cm} \caption{ $11$ crossing knots with a minimal quasi-alternating diagram.} \label{T22}
\end{table}

\begin{figure}
\begin{center}
\scalebox{.3}{
\includegraphics{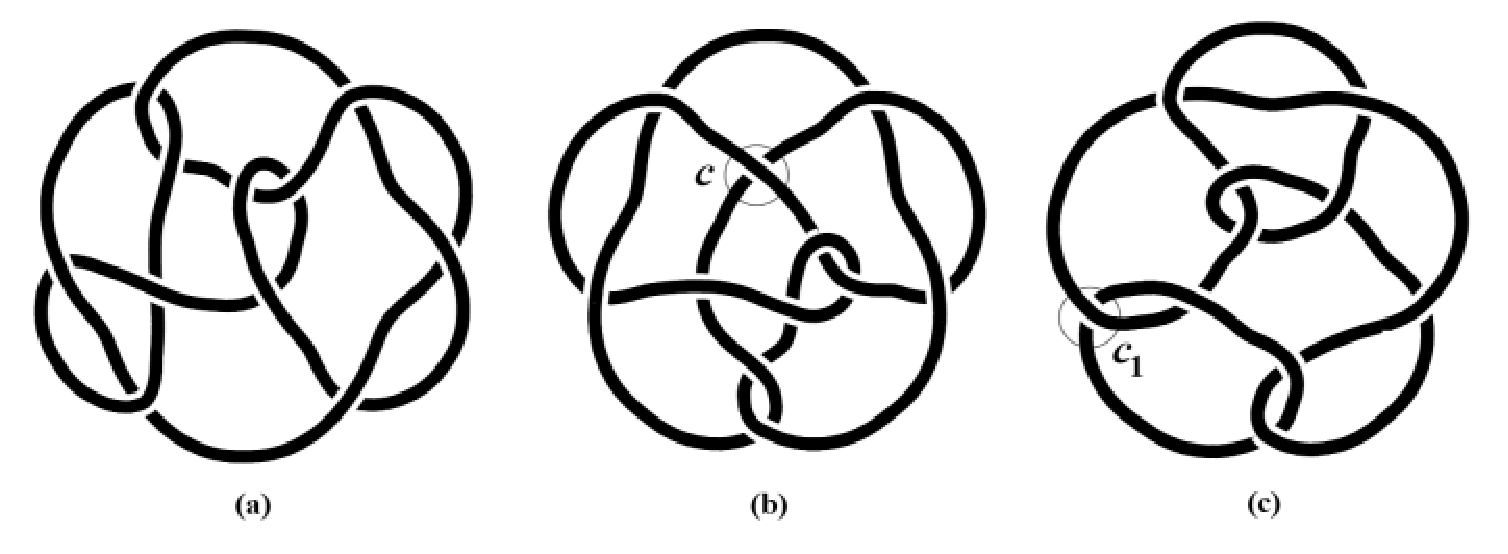}}
\caption{(a) Non-QA minimal diagram $(3,-2\,1)\,(2\,1,2)$ of the
knot $11n_{65}$; (b) the other QA minimal diagram
$6^*2.2\,1.-2\,0.-1.-2$ of the knot $11n_{65}$; (c) QA link
$(2,2+)\,-(2\,1,2)$.}  \label{L11n65}
\end{center} \end{figure}

The knot $11n_{65}$ has two minimal diagrams: $(3,-2\,1)\,(2\,1,2)$
(Fig. \ref{L11n65}a) and $6^*2.2\,1.-2\,0.-1.-2$ (Fig.
\ref{L11n65}b). The first diagram is not QA, and the second is QA.
By smoothing at the crossing $c$, the second diagram resolves into
unknot and QA link $(2,2+)\,-(2\,1,2)$, which resolves into QA knots
$3,2\,1,-2$ and $2\,1\,1,2\,1,-2$ by smoothing the crossing $c_1$
(Fig. \ref{L11n65}c). Moreover, all minimal $KL$ diagrams of the
family derived from knot $11n_{65}$, $(3,-2\,1)\,(p\,1,2)$ and
$6^*2.p\,1.-2\,0.-1.-2$ ($p\ge 2$) which represent the same $KL$
have this property: the first is not QA, and the other is QA.

The remaining Montesinos knot

$$11n_{50}=-2\,2,2\,2,3=M(0;(5,-2),(5,2),(3,1))=M(1;(5,3),(5,2),(3,1))$$

\begin{figure}
\begin{center}
\scalebox{.4}{
\includegraphics{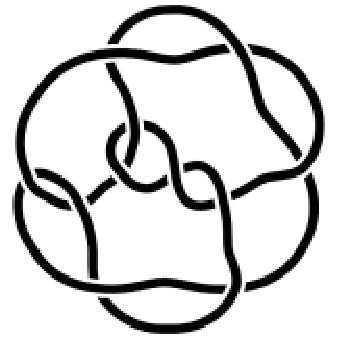}
} \caption{Knot $11n_{50}$. \label{L11n50}}
 \end{center} \end{figure}

\noindent (Fig.\ref{L11n50}) is homologically thin, and with no
minimal quasi-alternating diagrams.  In the first version of this
paper we proposed this knot as the smallest candidate for an
homologically thin knot which is not quasi-alternating. J. Greene
[Gr] proved that $11n_{50}$ is not QA, making it the  first example of
a homologically thin knot which is not QA.

$12$-crossing knots $12n_{196}$ = $(-3\,1,3)\,(2\,1,2)$, $12n_{393}$
= $8^*2.2\,0:-2\,1\,0$ and $12n_{397}$ = $2\,1\,1:-2\,1\,0:2\,0$
have another minimal diagram which is QA, and the knots given in
Table \ref{Lost} are candidates for homologically thin non-QA knots.

\begin{table}[h]
\begin{tabular}{|c|c|c|c|} \hline
 $12n_{139}$  & $.2.(-2\,1,2).2$ & $12n_{331}$  & $(-3,-2\,-1)\,(3,2+)$ \\ \hline
 $12n_{397}$  & $2\,1\,1:-2\,-1\,0:2\,0$ &   $12n_{414}$  &   $-2\,-1\,0.3.2.2\,0$      \\ \hline
 $12n_{768}$ & $2:-3\,-1\,0:3\,0$ & $12n_{838}$  & $-2.-2.-2\,0.2.2.2\,0$  \\ \hline
\end{tabular}\vspace{0.2cm}
\caption{Candidates for 12-crossing homologically thin non-QA
knots.} \label{Lost}
\end{table}

Kanenobu knots $K(n,3-n)$ are knots of the form
$2.2.-p\,0.-2.-2.-q\,0$, where $|p-q|=3$, $n\ge 0$. J.~Greene proved
that except knot $11n_{132}=K(1,2)=K(2,1)$ all Kanenobu knots
$K(n,3-n)$ are homologically thin non-QA knots. The first member of
this family is $11n_{50}=K(0,3)$, and the next members are knots
$12n_{414} = -2\,1\,0.3.2.2\,0$, $15_n{54616}$, {\it etc.}

All other 141 $\overline{Kh}'$-thin knots with $n=11$ crossings  and
656 $\overline{Kh}'$-thin knots with $n=12$ crossings crossings are
QA.

As the additional candidates for non-QA homologically thin knots, we
propose the following 13-crossing Montesinos knots, but we have only check their minimal diagram stated by the Conway symbol:

\begin{itemize}
\item $13n_{1408}=-3\,2,2\,2,2\,2=M(0;(7,-3),(5,2),(5,2))=M(1;(7,4),(5,2),(5,2))$,
\item $13n_{2006}=-3\,2,3\,2,3 =M(0;(7,-3),(7,3),(3,1))=M(1;(7,4),(7,3),(3,1))$, and
\item $13n_{3142}=-2\,2,3\,1\,2,3=M(0;(5,-2),(11,4),(3,1))=M(1;(5,3),(11,4),(3,1))$.
\end{itemize}

\section{Quasi-alternating links up to $12$ crossings}

In this section we study links up to $12$ crossings.  First we provide two examples of QA links with two distinct minimal diagrams,
where the first is not QA, and the other is: $(3,-2\,1)\,(2,2)$ =
$.2.-3\,0.2$, $(-2\,1,4)\,(2,2)$ = $6^*-3.-2.2\,0:2.-1$, {\it etc.}

Next we provide links with $11$ and $12$ crossings with no minimal QA diagrams.  Among $11$-crossing links there are five candidates for
homologically thin non QA links: $L11n_{77}=6^*2.(2,-2):2\,0$,
$L11n_{226}=6^*-2.2.-2:2\,1$, $L11n_{432}=8^*-2\,0:-2\,0:.2\,0$, and
$L11n_{437}=6^*-2.2.-2:2\,1\,0$. There are  eight candidates that are 12-crossing links:

\noindent $2\,1\,1,2\,2,-2\,1\,1$, $(-2\,1,2\,1\,1)\,(2,2+)$, $(-2
\,1,-4)\,(2,2+)$, $(2\,1,3),2,(2,-2)$, $6^*2:.(-2\,1,3)\,0$,
$6^*3.(2,-2):2\,0$, $8^*2.-2\,1\,0::2$, $8^*(2,-2)::-2\,0$.

J.~Greene [Gr] proved that link $L11n_{90}=(2,2+)\,-(2\,1,3)$ is
homologically-thin non QA link.

Moreover, the following families of homologically thin
links\footnote{Odd homology was computed for these families up to 16
crossings.} are potentially infinite families of non QA
homologically thin $KL$s: $(p,2+)\,-(2\,1,3)$ ($p\ge 2$),
$6^*p.(2,-2):2\,0$ ($p\ge 2$). Members of the first family
(beginning with the 11-crossing link $(2,2+)\,-(2\,1,3)$ and knot
$12n_{331}$  $=$ $(3,2+)\,-(2\,1,3)$) are links for even $p$ and
knots for odd $p$, and the members of the other family are
2-component links.

\section{Families of odd-homology thick links}

According to M.~Khovanov \cite{Kh2}, there is no doubt that for
large $n$, most $n$-crossing links are $\overline{Kh}$-thick. In the
case of odd Khovanov homology $\overline{Kh}'$, preliminary
computational results show that some families of thin links can
become thick.  The family of
$\overline{Kh}'$-thin knots $-2\,2,2\,2,p$ ($2\le p\le 5$) is a typical example- these
$\overline{Kh}'$-thick knots have torsion of order 5 for $p\ge 6$.

\begin{table}[h]
\begin{tabular}{|c|c|c|c|c|c|}   \hline
   No. of crossings & 6 & 7 & 8 & 9 & 10  \\ \hline
   No. of $\overline{Kh}'$-thick $KL$s & 1 & 1 & 6 & 15 & 61  \\ \hline
   No. of links which generate & & & & & \\
   new families of $\overline{Kh}'$-thick links & 1 & 0 & 4 & 6 & 28  \\ \hline
\end{tabular}
\medskip \caption{The number of $\overline{Kh}'$-thick links for
$6\le n\le 10$ and numbers of links which generate new families of $\overline{Kh}'$-thick links.} \label{TT4}
\end{table}

According to the computational results, most of the
$\overline{Kh}'$-thick links are members of the families of
$\overline{Kh}'$-thick links obtained for some lower value of $n$.
Table \ref{TT4} contains the number of $\overline{Kh}'$-thick links for
$6\le n\le 10$, and the
number of links which generate new families of $\overline{Kh}'$-thick links.\\

\begin{table}[h]
\begin{tabular}{|c|c|c|c|}  \hline
     $p\,q,r,s-$ & $p,q,r,s-$ &
     $(p,q)\,-(r,s)$ & $p,q,r,s--$  \\ \hline
\end{tabular}
\vspace{0.2cm} \caption{Families of $\overline{Kh}'$-thick links
for  minimal number of crossings $n \geq 8$ crossings.} \label{TT5}
\end{table}

Recall that families of $KL$s given in Conway notation
\cite{Con,Cau,JaSa} are determined by their generating links whose
Conway symbol contains integer parameters $p,q,r, \ldots$ greater or
equal to 2. There is only one family $p,q,r-$ of odd homology thick
links beginning for $n=6$ crossings. Tables
~\ref{TT5},~\ref{TT6},~\ref{TT7} contain families of
$\overline{Kh}'$-thick knots beginning from $n=8,9,10$ crossings,
together with additional conditions\footnote{Conditions are based on
the computational results for $KL$s with at most $16$ crossing in
each of the families.} on parameters insuring that the family
contains only $\overline{Kh}'$-thick links.

\begin{table}[h]
\begin{tabular}{|c|c|}   \hline
   $p\,1,q,r-$; $min(q,r)>p$ & $p\,1\,q,r,s-$ \\ \hline
     $(p,q+)\,(r,s-)$; $max(p,q)\le min(r,s)$ & $(p\,1,q)\,-(r,s)$  \\ \hline
     $(p,q+)\,-(r,s)$ & $.-(p,q)$  \\ \hline
\end{tabular}
\vspace{0.2cm} \caption{Families of $\overline{Kh}'$-thick links
beginning from $n=9$ crossings.} \label{TT6}
\end{table}

\begin{table}[b]
\begin{tabular}{|c|c|}   \hline
   $p\,q,r\,s,t-$ & $p\,1\,1\,q,r,s-$ \\ \hline
    $p\,1\,1,q,r-$; $min(q,r)\ge 3$ & $p\,1,q,r,s--$   \\ \hline
     $p,q,r,s,t-$ & $p,q,r,s,t--$ \\ \hline
     $p\,q,r,s,t-$ & $p\,q\,r,s,t-$  \\ \hline
     $p\,q,r,s,t--$ & $(p\,1,q)\,-(r\,1,s)$  \\ \hline
     $(p\,1\,1,q)\,-(r,s)$ & $(p,q,r--)\,(s,t)$  \\ \hline
     $(p\,q,r)\,-(s,t)$ & $(p,q)\,r\,(s,t-)$  \\ \hline
     $(p,q,r)\,(s,t-)$; $max(p,q,r)\le min(s,t)$ & $(p,q,r)\,-(s,t)$  \\ \hline
     $(p,q+r)\,(s,t-)$ &  $(p\,1,q+)\,(r,s-)$ \\ \hline
     $(p,q)\,-r\,(s,t)$ & $(p,q),r,-(s,t)$  \\ \hline
     $(p,q),r,(s,t-)$ & $p:-q\,0:-r\,0$; $min(p,q,r)\ge 3$  \\ \hline
     $-p\,1\,0:q\,0:r\,0$ & $.-(p,q).r$  \\ \hline
     $.-(p,q).r\,0$ & $.-(p,q):r\,0$  \\ \hline
     $.-(p,q):r$ & $.(p,q-)\,1$  \\ \hline
\end{tabular}
\vspace{0.2cm} \caption{Families of $\overline{Kh}'$-thick links
beginning from $n=10$ crossings.} \label{TT7}
\end{table}

 Inspired by the work of M. ~Khovanov on Ap-special knots \cite{Kh2}, as analogy we introduce Jp-special $KL$s.
 The Jones polynomial $J(L)=\sum c_iq^{2i}$ of a link $L$ is
\emph{alternating}, if $c_ic_j>0$ means that $j=i \mod 2$ and if
$c_ic_j<0$, then $j\neq i \mod 2$. The Jones polynomial has no gaps
if $c_i\neq 0$, $c_{i+k}\neq 0$ implies $c_{i+m}\neq 0$ for all $m$
between 1 and $k-1$. A link is called {\it Jp-special} if its Jones
polynomial is either non-alternating or has gaps.

\begin{table}[h]
\centering
\begin{tabular}{|c|c|}   \hline
   $(p,q)\,-(r,s)$  & $p,q,r,s--$ \\ \hline
   $(p\,1,q)\,-(r,s)$  & $.-(p,q)$\\ \hline
   $p\,1,q,r,s--$ & $p,q,r,s,t--$  \\ \hline
   $p\,q,r,s,t--$ & $(p\,1,q)\,-(r\,1,s)$\\ \hline
   $(p\,1\,1,q)\,-(r,s)$  & $(p,q,r--)\,(s,t)$ \\ \hline
   $(p\,q,r)\,-(s,t)$ &  $(p,q,r)\,-(s,t)$ \\ \hline
   $(p,q),r,-(s,t)$  & $.-(p,q).r$ \\ \hline
   $.-(p,q).r\,0$  & $.-(p,q):r\,0$ \\ \hline
   $.-(p,q):r$  &   \\ \hline
\end{tabular}
\vspace{0.2cm} \caption{Families of adequate Jp-special links.}
\label{T8}
\end{table}

Adequacy of links is determined using program {\it LinKnot}. Among
all above-listed families of $\overline{Kh}'$-thick links, families
that consist exclusively  of Jp-special links are the families of
adequate links given in Table \ref{T8}, and three families of
semi-adequate Jp-special links:

\begin{itemize}
\item $p,q,r-$; \item $p:-q\,0:-r\,0$ with $min(p,q,r)\ge 3$; and
\item $-p\,1\,0:q\,0:r\,0$.
\end{itemize}

\begin{con} All alternating knots or links, except those belonging to the family given by Conway symbol
$n$, $($i.e. $2_1^2$, $3_1$ ,$4_1^2$, $\ldots $ consisting of
torus links $(2,n)$, $n\ge 2$$)$, are not Jp-special.
\end{con}

\begin{table}[h]
\centering
\begin{tabular}{|c|c|}   \hline
   $p\,1,q\,1,r,s--$  & $p\,1,q,r\,1,s--$ \\ \hline
   $p\,1,q\,1,r\,1,s--$  & $(p\,1\,1,q) -(r\,1,s)$ \\ \hline
   $(p\,q,r) -(s\,1,t)$  & $.-(p\,1,q).r$ \\ \hline
   $.p.-(q\,1,r)$  & $6^*-(p\,1,q).r\,0$ \\ \hline
   $6^*-(p\,1,q).r\,0$  & $6^*p.-(q,r\,1)$ \\ \hline
   $p\,1,q,r,s,t--$  & $p\,1\,q,r,s,t--$ \\ \hline
   $p\,1\,1,q,r,s--$  & $(p\,1\,1\,1,q)\,-(r,s)$ \\ \hline
   $(p\,1\,q,r)\,-(s,t)$  & $(p\,q\,1,r)\,-(s,t)$ \\ \hline
   $(p,q\,1)\,-(r\,1,s\,1)$  & $(p\,1,q\,1\,1)\,-(r,s)$ \\ \hline
   $(p\,1,q\,r)\,-(s,t)$  & $(p\,1,q\,1)\,-(r,s)$ \\ \hline
   $(p,q)\,-(r\,1,s\,1)$  & $(p\,1,q,r--)\,(s,t)$ \\ \hline
   $(p\,1,q,r)\,-(s,t)$  & $(p,q,r--)\,(s,t)$ \\ \hline
   $(p\,1,q,r--)\,(s,t)$  & $(p,q\,1,r--)\,(s,t)$ \\ \hline
   $(p,q\,1,r)\,-(s,t)$  & $(p,q\,1,r--)\,(s,t)$ \\ \hline
   $(p,q,r)\,-(s\,1,t)$  & $(p,q,r--)\,(s\,1,t)$ \\ \hline
   $(p,q,r--)\,(s,t+)$  & $-(p\,1,q),r,(s,t)$ \\ \hline
   $(p\,1,q),r,-(s,t)$  & $(p,q),r\,1,-(s,t)$ \\ \hline
   $-(p,q),r,(s,t+)$  & $(p,q)\,1,-(r,s),t$ \\ \hline
   $6^*-(p\,1,q\,1)$  & $6^*-(p,q\,1\,1)$ \\ \hline
   $6^*-(p,q r)$  & $6^*-(p,q\,1)$ \\ \hline
   $6^*-(p\,1,q)$  & $6^*-(p,q),-r$ \\ \hline
   $6^*-(p,q),r$  & $6^*(p,q,r--)$ \\ \hline
   $6^*-(p,q,r)$  & $6^*p\,1:.-(q,r)\,0$ \\ \hline
   $6^*p:.-(q,r)\,0$  & $6^*2:.-(p\,1,q)\,0$ \\ \hline
   $6^*p:.-(q,r)\,0$  & $6^*p\,1:.-(q,r)$ \\ \hline
   $6^*p:.-(q,r)$  & $6^*p:.-(q\,1,r)$ \\ \hline
   $6^*p:.-(q,r)$  & $6^*-(p,q).r\,0.s$ \\ \hline
   $6^*-(p,q).r.s\,0$  & $6^*p.-(q,r).s\,0$ \\ \hline
   $6^*-(p,q).r\,0::s\,0$  & $6^*-(p,q).r.s$ \\ \hline
   $6^*p.-(q,r).s$  & $8^*-(p,q)\,0$ \\ \hline
   $8^*-(p,q)$  & $10^{***}p::-1.-1.-1.-1:-1$ \\ \hline
   $10^{***}::-1.-1.-1.-1.p\,0.-1$  & $$ \\ \hline
\end{tabular}
\vspace{0.3cm}
 \caption{Familes of adequate $\overline{Kh}'$-thick
links beginning from $n=11$ crossings.} \label{T9}
\end{table}

\begin{con} Every adequate non-alternating link is Jp-special.
\end{con}

\begin{con}  All minimal positive braids are Jp-special.
\end{con}

Among links with $11$ crossings at least $126$ generate families of
$\overline{Kh}'$-thick links. Sixty-one of these families consist
only of adequate links, see Table \ref{T9}.

The remaining 65 families containing $\overline{Kh}'$-thick links
are given in the following table, together with estimated conditions
for $\overline{Kh}'$-thick links, see Table \ref{T6}.

\begin{table}[h]
\centering
\begin{tabular}{|c|c|}   \hline
   $p\,1\,1\,1,q,r-;$ $min(q,r)\ge 3$  & $p\,1\,1\,q,r,s$ \\ \hline
   $-p\,1\,1:q:r;$ $p\le min(q,r)$ or $min(q,r)\ge 3$ & $-p\,1\,1\,0:q\,0:r\,0$ \\ \hline
   $p.q\,0.-r.s.t\,0;$ $min(p,q)\ge r$ & $-p\,0:q\,1:-r\,0$\\ \hline
   $p\,q:-r\,0:-s\,0$  & $p.-q\,1.-r\,0.s\,0$ \\ \hline
   $p.-q.-r\,0.s\,0$  & $8^*p.-q\,0.r$ \\ \hline
   $p\,q,r\,1,s-;$ $q>r$ and $s>r$ & $9^*.-p:.-q$ \\ \hline
   $p\,q\,1 r,s,t-$  & $p\,1\,q\,r,s,t-$ \\ \hline
   $p\,1\,q,r\,s,t-$  & $p\,1\,1\,1\,q,r,s-$ \\ \hline
   $(p,q-)\,(r\,1,s\,1+)$   & $-(p\,q,r)\,(s,t+)$\\ \hline
   $(p\,q,r-)\,(s,t+);$ $r\ge t$ & $-(p,q)\,(r\,s,t+)$ \\ \hline
   $(p,q-)\,(r\,1\,1,s+);$ $p\ge s$  & $(p,q-)\,(r\,1,s+t)$ \\ \hline
   $(p\,1,q-)\,-1\,-1\,(r,s)$ & $(p\,1,q)\,-1\,-1\,(r,s-)$ \\ \hline
   $(p\,1,q)\,1\,1\,-(r,s)$  & $-(p,q)\,1\,1\,(r,s)$ \\ \hline
   $(p,q)\,1\,1\,-(r,s)$ &  $-(p,q)\,1\,1\,(r,s+)$ \\ \hline
   $(p\,1,q,r)\,(s,t-);$ $max(q,r)\le min(s,t)$ & $(p,q,r+)\,(s,t-)$ \\ \hline
   $(p,q,r+)\,-(s,t);$ $min(p,q,r)\ge max(s,t)$ & $(p,q,r-)\,(s,t+)$  \\ \hline
   $-(p,q,r)\,(s,t+)$  & $(p,q),r\,1,(s,t-)$ \\ \hline
   $(p,q-),r+,(s,t-)$   & $(p,q),r+,(s,t-)$ \\ \hline
   $(p,q),-r,-1,(s,t)$ & $(p,q)\,-1\,-1\,-1\,(r,s)$ \\  \hline
   $(p,q-)\,1\,1\,1\,(r,s-)$  & $(p,q)\,-1\,-1\,-1\,(r,s-)$ \\ \hline
   $(p,q)\,1\,1\,1\,-(r,s)$  &  $(p,q)\,r\,1\,-(s,t)$ \\ \hline
   $6^*-p\,q.r\,0.s;$ $q\le min(r,s)$ &  $6^*-p\,1.q.-r$ \\ \hline
   $6^*-p.q\,0.r;$ $q\ge p$ and $r\ge p$ & $6^*-p.q\,1.-r\,1$ \\ \hline
   $6^*-p.q.-r;$ $min(p,q,r)>2$ & $6^*p.q.r\,0.-s\,1$ \\ \hline
   $6^*p.q.-r\,1.s\,0$  & $6^*-p\,1.q.-r:s$ \\ \hline
   $6^*p.q.-r.s\,0.t;$ $min(p,q)\ge r$ and $min(p,q,r)\le min(s,t)$ & $6^*-p\,1.q\,0.r.s\,0$ \\ \hline
   $6^*-p.q.-r:s;$ $min(p,q,r)\ge 3$ &  $6^*p.q\,0.-r.s\,0.t$ \\ \hline
   $6^*(p,q-);$ $min(p,q)\ge 3$ &  $6^*(p,q-)\,r$  \\ \hline
   $6^*(p,q-),r$  & $6^*p.-(q,r)\,1$ \\ \hline
   $6^*p:.-(q,r)\,1\,0$  & $6^*p:.(q,r-)\,1\,0$  \\ \hline
   $6^*(p,q).-r\,0.-s$ & $6^*-(p,q).r.-s$ \\ \hline
   $6^*p.(q,r-).s$  & $8^*p.-q\,1\,0$ \\ \hline
   $8^*-p\,0.-q\,0.-r\,0$  &  \\ \hline
\end{tabular}
\vspace{0.5cm} \caption{The remaining 65 families containing
$\overline{Kh}'$-thick links.} \label{T6}
\end{table}

$19$ families in Table \ref{T6}  contain only Jp-special
links, see Table \ref{T7}.
\begin{table}[h]
\centering
\begin{tabular}{|c|c|}   \hline
   $-p\,1\,1\,0:q\,0:r\,0$  & $-p\,0:q\,1:-r\,0$ \\ \hline
   $p\,q:-r\,0:-s\,0$  & $p.-q.-r\,0.s\,0$ \\ \hline
   $9^*.-p:.-q$  & $(p,q-)\,-1\,-1\,(r,s)$ \\ \hline
   $(p\,1,q)\,-1\,-1\,(r,s-)$  & $(p,q-),r+,(s,t-)$ \\ \hline
   $(p,q),-r,-1,(s,t)$  & $(p,q)\,-1\,-1\,-1\,(r,s)$ \\ \hline
   $6^*-p\,1.q.-r$  & $6^*-p.q\,1.-r\,1$ \\ \hline
   $6^*-p.q.-r$  & $6^*p.q.r\,0.-s\,1$ \\ \hline
   $6^*p.q.-r\,1.s\,0$  & $6^*-p\,1.q\,0.r.s\,0$ \\ \hline
   $6^*(p,q).-r\,0.-s$  & $6^*-(p,q).r.-s$ \\ \hline
   $8^*-p\,0.-q\,0.-r\,0$  & $$ \\ \hline
\end{tabular}
\vspace{0.5cm} \caption{Families of Jp-special links from the
Table 6} \label{T7}
\end{table}

The majority of the remaining families derived from
$\overline{Kh}'$-thick links with $n=11$ crossings contain both
Jp-special and not Jp-special links, but two families, $6^*(p,q-)$
($p\ge 3$, $q\ge 3$) and $8^*p.-q\,1\,0$, contain only links which
are not Jp-special.


\section{Recognition of odd-homology thickness  based on
computational results for different classes of links}

After computing odd homology of different classes of links by
using the program {\it KhoHo} \cite{Sh1}, we propose several
conjectures about odd homology thick links.

\begin{con} Every link given by a positive minimal
$k$-braid is $\overline{Kh}'$-thick $(k\ge 3)$\footnote{See
\cite{Kh2}, Problem 6.2.}. \label{MinBraid1}
\end{con}

According to computational results, this conjecture holds for all
positive minimal $3$-, $4$- and $5$-braids with at most $n=20$
crossings.

If Conjecture \ref{MinBraid1} is true in general, it implies that all
non-alternating Lorenz links \cite{Lo,BiWi,Gh,GhLe} are
$\overline{Kh}'$-thick{\footnote{Only alternating Lorenz links are
obtained from minimal $2$-braids of the form $a^n$, giving the
family of the links $2_1^1$, $3_1$, $4_1^2$, $5_1$,...}. M.~Sto\v
si\'c proved that non-alternating torus knots, a subset of Lorenz
links, are $\overline{Kh}$-thick \cite{St}.

A rational tangle is called {\it positive} if it has only positive
numbers in its Conway symbol. Recall that up to taking mirror image
and permuting rational tangles, every non-alternating Montesinos
link can be denoted by Conway symbol of the form
$p_1,...,p_m,-q_1,...,-q_n$ where $p_i$ and $q_j$ denote positive
rational tangles which do not start by one, $i=\in \{1,2,...,m\}$,
$j \in \{1,2,...,n\}$, $m\ge n\ge 1$, $m\ge 2$.

Every Montesinos link with $n>1$ is $\overline{Kh}'$-thick, since it
is an adequate non-alternating link \cite{Ja}. Hence, we restrict
our consideration to the Montesinos links of the form
$p_1,...,p_m,-q$. Every Montesinos link can be given by one or
several different Conway symbols with the minimal number of
crossings. For example, the link $2\,1\,2,2\,1\,1,-2\,1$ can be
written also as $2\,1\,2,-2\,2,3$ or $-2\,1\,1\,1,2\,1\,1,3$. A
rational tangle $t$ is of length $l$ if its symbol consists of $l$
integers, i.e., $t=t_1\,t_2\,...t_l$. To rational tangles of length
$l>1$ belonging to a Conway symbol of a Montesinos link we will
apply the following reduction: every positive rational tangle is
replaced by its last number $t_l$, every negative rational tangle by
$-t_l-1$, and tangles of the length $l=1$ remain unchanged. In this
way, from every Conway symbol of a Montesinos link $p_1,...,p_m,-q$
we obtain the reduced Conway symbol $\overline p_1,...,\overline
p_m,-\overline q$.

\begin{con}
Montesinos link of the form $p_1,...,p_m,-q$ is
$\overline{Kh}'$-thick if it has a Conway symbol with a minimal
number of crossings satisfying the relationship $min(\overline
p_1,\overline p_2,...,\overline p_m)\ge \overline q.$ \label{Con}
\end{con}

For pretzel links all tangles are of the length 1, so $\overline
p_i=p_i$ ($i=1,\ldots,m$), $\overline q=q$, and Conjecture 5.2
corresponds to Proposition 2.2 [Gr]:

\begin{proposition}
For $n\ge 2$ and $p_1$,$\ldots $,$p_n\ge 2$, and $q\ge 1$, the
pretzel link $P(p_1,\ldots ,p_n,-q)$ is QA iff $q>min\{p_1,\ldots
,p_n\}$ $[Gr]$.
\end{proposition}

For $n\le 11$ crossings all non-alternating thick Montesinos links
are exactly those described in Conjecture \ref{Con}. However, for
$n\ge 12$ crossings some $\overline{Kh}'$-thick links do not satisfy
the conditions of Conjecture \ref{Con}. For example, for $12\le n
\le 15$ crossings, we have the following families of exceptional
Montesinos links:

\begin{enumerate}
 \item $-2\,2,2\,2,p$; $p\ge 6$; \item $-3\,3,2\,3,p$; $p\ge 4$;
\item $-p\,2,q\,2,r$; $p>q\ge 2$, $r\ge 3$; \item
$-p\,2,q\,1\,2,r$; $p=2$, $q\ge 2$, $r\ge 4$ or $p>2$, $q\ge 2$,
$r\ge 3$, {\it etc.}
\end{enumerate}

Let $p_1,p_2,...,p_l-k$ ($l\ge 2$, $k\le l$) denote tangle
$p_1,p_2,...,p_l-...-$ with $k$ minuses, where $p_i$ ($i=1,2,...,l$)
are positive rational tangles that do not start with 1.

\begin{theorem}
All algebraic links of the form $(p_1,...,p_m)\,(q_1,...,q_n-k)$,
where $m,n \ge 2$ and all $p_i$, $i\in\{1,...,m\}$,  and $q_j$,
$j\in\{1,...,n\}$,  denote positive alternating rational tangles,
are $\overline{Kh}'$-thick for $k\ge 2$. \label{SpecCase}
\end{theorem}

The Theorem \ref{SpecCase} holds since all algebraic links above
are non-alternating adequate links \cite{Ja}, which are
homologically thick.

Hence, we can restrict our consideration to links of the form
$(p_1,...,p_m)\,(q_1,...,q_n-)$. First we consider the case where
$p_i$ and $q_j$ and are positive integer tangles different from 1
($i=1,2,...,m$; $j=1,2,...,n$; $m\ge 2$, $n\ge 2$).

\begin{con} A link of the form
$(p_1,...,p_m)\,(q_1,...,q_n-)$, where $p_i$ and $q_j$ are positive
integer tangles different from 1 $(i=1,2,...m$; $j=1,2,...,n$; $m\ge
2$, $n\ge 2)$ is $\overline{Kh}'$-thick if $m\ge 4$ or $m=3$ and
$max(p_1,...,p_m)\le min(q_1,...,q_n)$.
\end{con}

\begin{con}
A link of the form $(p_1,...,p_m)\,(q_1,...,q_n-)$, where $p_i$ and
$q_j$ are positive rational tangles that do not start with $1$
$(i=1,2,...m$; $j=1,2,...,n)$ which are not all integer tangles, is
$\overline{Kh}'$-thick if $m\ge 3$, $q=min(\overline q_1,\overline
q_2,...,\overline q_n)>1$ and
\begin{itemize}
\item if length$(p_i)\ge 2$ then $\overline p_i=1$
$(i=1,2,...,m)$, and \item $max(\overline p_1,\overline
p_2,...,\overline p_m)\le q$.
\end{itemize}
\end{con}

Next, we consider algebraic links of the form $(p_1,...,p_m+k)
\,(q_1,...,q_n-l)$, where $p_i$ ($i=1,...,m$, $m\ge 2$), $q_j$
($j=1,...,n$) are positive alternating rational tangles, $+k$
denotes a sequence of $k$ pluses ($k\ge 1$), and $-l$ denotes
sequence of $l$ minuses ($n>l\ge 1$). All these links are
$\overline{Kh}'$-thick for $l>1$, since they are non-alternating
adequate links. Hence, we need to analyze only links of the form
$(p_1,...,p_m+k)\,(q_1,...,q_n-)$. We propose the following
conjecture:

\begin{con} A link of the form
$(p_1,...,p_m+k)\,(q_1,...,q_n-)$ is $\overline{Kh}'$-thick if $m\ge
3$, $q=min(\overline q_1,\overline q_2,...,\overline q_n)>1$ and
\begin{itemize}
\item if $length(p_i)\ge 2$ then $\overline p_i=1$
$(i=1,2,...,m)$, and \item $max(\overline p_1,\overline
p_2,...,\overline p_m)+k\le q.$ \label{Con6}
\end{itemize}
\end{con}

\begin{con} An algebraic link of the form $(p_1,...,p_m+k)
-(q_1,...,q_n)$, where $p_i$ $(i=1,...,m$, $m\ge 2)$ and $q_j$
$(j=1,...,n$, $n\ge 2)$ are positive alternating rational tangles
and $+k$ denotes a sequence of $k$ pluses $(k\ge 1)$ is
$\overline{Kh}'$-thick if $q=min(\overline q_1,\overline
q_2,...,\overline q_n)>1$ and

\begin{itemize}
\item $min(\overline p_1,\overline p_2,...,\overline p_m)>1$,
$k<q$ or \item $min(\overline p_1,\overline p_2,...,\overline
p_m)=1$, $k+1<q.$ \label{Con7}
\end{itemize}
\end{con}

However, there are $\overline{Kh}'$-thick links belonging to this
class, which do not satisfy the conditions of Conjecture \ref{Con7},
for example links of the family $(2,2,2+)\,-(2\,1\,1,p)$ for $p\ge
3$.

Let $p_1,p_2,...,p_m-k$ ($m\ge 2$, $m\ge k$) denote tangle
$p_1,p_2,...,p_m-...-$ with $k$ minuses, where $p_i$ ($i=1,2,...,m$)
are positive rational tangles that do not start with 1.

\begin{theorem}
Every polyhedral link containing a tangle $p_1,p_2,...,p_m-k$ $(m\ge
2$, $m\ge k)$ in one vertex of the basic polyhedron, and arbitrary
positive tangles in other vertices is $\overline{Kh}'$-thick for
$k\ge 2$, since it is a non-alternating adequate link \cite{Ja}.
\end{theorem}

\bigskip

\noindent {\bf Acknowledgments} We are thankful to Alexander
Shumakovitch for providing us the new version of his remarkable
knot-theory computer program {\it KhoHo}, Alexander Stoimenow for
the list of non-minimal diagrams of the knot $11n_{50}$, and Etienne
Ghys and Pierre Dehornoy for the list of braids corresponding to
Lorenz knots and links.

\bigskip


\bigskip
\bigskip

\footnotesize

\noindent THE MATHEMATICAL INSTITUTE, KNEZ MIHAILOVA  36, P.O.BOX
367, \\ 11001 BELGRADE, SERBIA

\medskip

\noindent {\it E-mail address:} $\mathrm{sjablan@gmail.com}$

\bigskip

\noindent THE MATHEMATICAL SCIENCES RESEARCH INSTITUTE, 17 GAUSS
WAY, BERKELEY, CA 94720-5070, USA

\medskip

\noindent {\it E-mail address:} $\mathrm{rsazdanovic@math.ncsu.edu}$

\end{document}